\documentclass{article}

\begin{document}

\title{ Abel Prize 2009 as a necessitation to discussion on geometry, or
physicist versus geometers.}
\author{Yuri A. Rylov}
\date{Institute for Problems in Mechanics, Russian Academy of Sciences,\\
101-1, Vernadskii Ave., Moscow, 119526, Russia.\\
e-mail: rylov@ipmnet.ru\\
Web site: {$http://rsfq1.physics.sunysb.edu/\symbol{126}rylov/yrylov.htm$}\\
or mirror Web site: {$http://gasdyn-ipm.ipmnet.ru/\symbol{126}%
rylov/yrylov.htm$}}
\maketitle

\begin{abstract}
Contemporary geometers do not acknowledge nonaxomatizable
geometries. This fact means that our knowledge of geometry is
poor. A perfect knowledge of geometry is important for "consumers
of geometry"\ (physicists dealing with geometry of microcosm),
because poor knowledge of geometry forces physicists to refuse
from the orthodox geometrical paradigm and to use the quantum
paradigm, which is a simple fitting, compensating our imperfect
knowledge of geometry.
\end{abstract}

Situation in geometry, when we know only a negligible part of geometries
serviceable for description of the space-time, should be qualified as a
\textit{poor knowledge of geometry}. Any geometry is a continual set of
propositions on properties of geometrical objects. A geometry is
axiomatizable, if the whole continual set of proposition can be deduced from
a finite set of basic propositions (axioms) by means of rules of the formal
logic. We know only axiomatizable geometries, because we do not enable to
construct nonaxiomatizable geometries. However, our inaptitude of
nonaxiomatizable geometries construction does not mean, that
nonaxiomatizable geometries do not exist.

Of late one suggested a method of construction of physical geometries, that
is geometries, which are described in terms and only in terms of the world
function \cite{R2002,R2007,R2008}. (The world function $\sigma \left(
P,Q\right) =\frac{1}{2}\rho ^{2}\left( P,Q\right) $, where $\rho \left(
P,Q\right) $ is a distance between the points $P$ and $Q$). Such an approach
to geometry, when a geometry is described completely by a distance, is
correlated better with our apprehension that a geometry is a science on
disposition of geometrical objects, than with the conventional apprehension,
when geometry is simply a logical construction.

The physical geometries are nonaxiomatizable as a rule, because the
equivalence relation is intransitive in the physical geometries, in general.
The equivalence relation is always transitive in axiomatizable geometries,
and geometries with intransitive equivalence relation cannot be
axiomatizable.

Almost all contemporary geometers do not acknowledge physical geometries. As
far as I understand, the reason of this fact is very simple. A construction
of a physical geometry does not refer to means of the geometry construction,
which are characteristic for axiomatizable geometry (manifold, linear vector
space, coordinate system, dimension, continuity). A physical geometry is
formulated in terms of points and of distance (world function) between them.
The physical geometry is formulated uniformly on arbitrary set of points
(continuous and discrete, with unlimited divisibility and with restricted
divisibility). A use of the physical geometry is reasonable and natural, if
the geometry is considered as a science on disposition of geometrical
objects. However, contemporary geometers consider means of the geometry
description as necessary attributes of a geometry. They cannot imagine a
geometry without such an attribute as the linear vector space. They consider
only geometries given on a manifold and they cannot imagine a geometry,
given on a set of points without a mention of the set continuity.

In the nineteenth century all mathematicians did not acknowledge the
non-Euclidean geometry, suggested by Lobachevsky and Bolyai. What was a
reason of such abruption? I did not find any investigation of this question
in literature. I think as follows. In that time the geometers dealt only
with the Euclidean geometry. The Cartesian coordinate system can be
introduced in the Euclidean geometry, but it cannot be introduced in the
non-Euclidean geometry. The geometers considered the Cartesian coordinate
system as an attribute of any geometry, because it was an attribute of the
Euclidean geometry. The geometers did not acknowledge the non-Euclidean
geometry as a real geometry, because one cannot introduce Cartesian
coordinate system in it. It was the first crisis in geometry.

Now we have the second crisis. The contemporary geometers do not acknowledge
the physical (nonaxiomatizable) geometries. The reason of such a relation
is, in particular, the fact, that one cannot introduce a linear vector space
in a physical geometry. It is supposed conventionally, that the linear
vector space is a necessary attribute of any geometry. Accordingly to this
approach the contemporary geometers suppose that any geometry is
axiomatizable and try to construct nonaxiomatizable geometries, deducing
them from some axiomatics. As far as it is impossible, the geometrical
constructions appear to be inconsistent or overdetermined. For instance, the
Riemannian geometry appears to be inconsistent (overdetermined) \cite%
{R2005,R2006}.

Physicists, who use a geometry for description of the space-time, meet
serious problems, connected with poor knowledge of geometry. Indeed, we know
only negligible part of possible geometries. The real space-time geometry
belongs to that part of geometries, which we do not know, and we are forced
to invent different hypotheses, compensating a use of wrong space-time
geometry. It happened as follows.

Considering physical phenomena in the microcosm, we may neglect influence of
the matter distribution on the space-time geometry. Then the space-time
geometry is to be uniform and isotropic. In the twentieth century the
Riemannian geometry was considered as the most general geometry serviceable
for description of the space-time. The geometry of Minkowski is the only
Riemannian geometry, which is uniform and isotropic. Investigations show,
that classical principles of dynamics together with the space-time geometry
of Minkowski cannot explain experimental data. As far as our knowledge of
geometry were poor, and we did not know uniform isotropic geometries other,
than the geometry of Minkowski, we cannot modify space-time geometry, to
agree predictions of the theory with experimental data. We were forced to
modify principles of dynamics. As a result the quantum paradigm (variation
of dynamic principles + fixed space-time geometry) has been existing for the
twentieth century. The geometric paradigm (classical principles of dynamics
+ variated space-time geometry) could not appear until appearance of other
(nonaxiomatizable) uniform isotropic space-time geometries.

The geometric paradigm is more reasonable and natural, than the quantum
paradigm. The physical geometries are formalized in the sense, that all
possible geometries are labelled by the world function. Giving the world
function, one determines the geometry. It is necessary only to choose a
proper world function. Principles of dynamics are not formalized (they are
not labelled), and it is very difficult to guess true principles of
dynamics. The quantum principles of dynamics work well for the
nonrelativistic physical phenomena. Investigating the relativistic physical
phenomena, we come to geometrical problems (strings, branes,
compactification) again, and our poor knowledge of geometry cannot help us
\cite{R2008a}. I have described the situation from the viewpoint of a
physicist, who is a "consumer of a geometry"

The situation looks differently from the viewpoint of a geometer, who is a
creator of a geometry. A geometer ignores the consume side of a geometry. He
perceives the Riemannian geometry and the physical geometry as different
matters, because they are constructed by different methods, although both
geometries are sciences on the disposition of geometric objects. At the same
time the geometer perceives the Euclidean geometry and the symplectic
geometry as similar matters, because both geometries are logical
constructions. He perceives the symplectic geometry as a geometry, although
it has no relation to the geometry as a science on the disposition of
geometric objects. He united Euclidean geometry and symplectic geometry by
the same term 'geometry'. The fact, that the symplectic geometry is logical
structure, as the Euclidean geometry, is more important for a geometer, than
the fact, that the symplectic geometry has no relation to a geometry as a
science on disposition of geometrical objects.

Situation may be illustrated by the legend on Alexander the Great, who hew
the Gordian Knot by his sword instead of to untangle it. The Gordian cart
was attached to the chancel by means of the Gordian Knot. From practical
viewpoint (separation of the cart from the chancel) the action of Alexander
is quite right. However, Alexander the Great did not manifest the art of
untangling of knots, and he could not to awarded as a "champion of
untangling of knots". If the separation of the cart from the chancel is
restricted by the condition, that the Gordian Knot must be untangled, the
separation becomes quite another problem, than that one, which had been
solved by Alexander the Great.

In a like way geometers consider the problem of a geometry construction only
under constraint, that this geometry is axiomatizable and it can be
considered to be a logical construction. It does not mean, that geometers do
not understand, that a geometry may be a science on disposition of
geometrical objects. Geometers understand that a geometry may be a science
on disposition of geometric objects, which is described in terms of only
distance. There exists so called "distance geometry" \cite{M28,B53}.
However, Blumental failed to construct pure metric geometry. He did not know
the deformation principle. He was forced to introduce additional concept of
a curve, which cannot be expressed via the concept of distance. As a result
now geometers consider any geometry as a logical construction. In other
words, any geometry is constructed under additional constraint, that the
nonaxiomatizable geometries are impossible.

The Abel Prize 2009 (on geometry) has been awarded under this unfounded
constraint, that the geometry may be only axiomatizable. Such approach would
be justified twenty years ago, when nonaxiomatizable geometries were not
known. Now, when such geometries appear to be possible \cite%
{R2002,R2007,R2008}, the Abel Prize to a geometer, dealing only with
axiomatizable geometries looks as a stimulation of the art of deduction and
of the art of complicated mathematical calculations. Of course, the art of
deduction should be stimulated, awarding prizes to mathematicians. However,
the awarder should declare, that the prize is awarded for the art of
mathematical calculations and for the art of deduction, but not for a
progress in the investigations of a geometry. The reason of awarding is to
be clear for all people. The most people perceive the geometry as "consumers
of a geometry", but not as its creators. From the viewpoint of the geometry
consumer (physicist) the real progress in the geometry investigation is
possible only, if the restriction on the geometry axiomatizability is
removed. Attempt of  the axiomatizable geometries construction, when even
the "axiomatizable" Riemannian geometry appears to be inconsistent, is not a
progress in geometry.

\end{document}